\documentclass[titlepage]{amas}
\tikzcdset{arrow style=Latin Modern}
\usepackage{amssymb, amsmath, amsfonts}
\usepackage[croatian]{babel}

\renewcommand{\dj}{d\kern-0.4em\char"16\kern-0.1em}

\usepackage{bm}
\def\genfd{{\bm k}}
\def\id{{\rm id}}
\long\def\nodo#1{{}}
\def\gg{\mathfrak{g}}
\def\bbC{\mathcal{C}}

\begin{document}
\title{A note on symmetric orderings}

\author{Zoran \v{S}koda}
\address{(1) Department of Teachers' Education, University of Zadar, Franje Tu\dj mana 24, 23000 Zadar, Croatia\\
(2) Division of Theoretical Physics, Ru\dj er Bo\v skovi\'c Institute, Bijeni\v cka cesta 54, 10002 Zagreb, Croatia}
\email{zskoda@unizd.hr}
\maketitle
\begin{abstract}
  Let $\hat{A}_n$ be the completion by the degree of a differential operator
  of the $n$-th Weyl algebra
  with generators  $x_1,\ldots,x_n,\partial^1,\ldots,\partial^n$.
  Consider $n$ elements $X_1,\ldots,X_n$ in $\hat{A}_n$ of the form
  $$
  X_i = x_i + \sum_{K = 1}^\infty \sum_{l = 1}^n\sum_{j = 1}^n x_l p_{ij}^{K-1,l}(\partial)\partial^j,
  $$
where $p^{K-1,l}_{ij}(\partial)$ is a degree $(K-1)$ homogeneous polynomial in $\partial^1,\ldots,\partial^n$, antisymmetric in subscripts $i,j$. Then for any natural $k$ and any function $i : \{1,\ldots,k\}\to\{1,\ldots,n\}$ we prove
  $$
\sum_{\sigma \in \Sigma(k)}
X_{i_{\sigma(1)}}\cdots X_{i_{\sigma(k)}}\triangleright 1 =  k! \,x_{i_1}\cdots x_{i_k},
$$
where $\Sigma(k)$ is the symmetric group on $k$ letters and
$\triangleright$ denotes the Fock action of the $\hat{A}_n$ on the
space of (commutative) polynomials. 
\end{abstract}

\section{Introduction and motivation}
In an earlier article~\cite{ldWeyl},
we derived a universal formula for an embedding of
the universal enveloping algebra $U(\mathfrak{g})$
of any Lie algebra $\mathfrak{g}$ with underlying rank $n$
free module over a commutative ring $\genfd$ containing the field $\mathbb{Q}$
of rational numbers into a completion $\hat{A}_{n,\genfd}$
of the $n$-th Weyl algebra over $\genfd$. 

\begin{definition}
The $n$-th Weyl algebra $A_{n,\genfd}$ over a commutative ring $\genfd$ is the associative $\genfd$-algebra defined by generators and relations as follows:
  $$
  A_{n,\genfd} := \genfd\langle x_1,\ldots,x_n,\partial^1,\ldots,\partial^n\rangle/
  \langle [x_i,x_j ],[\partial^i,\partial^j],[x_i,\partial^j]-\delta^j_i,\,i,j,=1,\ldots,n\rangle.
  $$
\end{definition}

We use the ``contravariant'' notation
for the generators of $A_{n,\genfd}$~(\cite{halg}, 1.1)
and $\delta^i_j$ is the Kronecker symbol.
The reader should recall the usual interpretation of the Weyl algebra elements as regular differential operators~\cite{coutinho,ldWeyl}. In other words,
the elements of $A_{n,\genfd}$ act on the polynomial algebra $\genfd[x_1,\ldots,x_n]$, consisting of commutative polynomials via 
the physicists' Fock action here denoted
by $\triangleright : A_{n,\genfd}
\otimes\genfd[x_1,\ldots,x_n]\to \genfd[x_1,\ldots,x_n]$. By definition, generators $x_i$ act as the multiplication operators by $x_i$ and
$\partial^j$ act as partial derivatives. 
The unit polynomial $1\in\genfd[x_1,\ldots,x_n]$ is interpreted as the vacuum state.

Complete $A_{n,\genfd}$ along the filtration given by the degree of differential operator~(\cite{ldWeyl,halg,leib}); the completion will be denoted $\hat{A}_{n,\genfd}$. Thus, the elements in $\hat{A}_{n,\genfd}$ can be represented as arbitrary power series in $\partial^1,\ldots,\partial^n$ with coefficients (say on the left) in the polynomial ring $\genfd[x_1,\ldots,x_n]$.

For a fixed basis $X_1^{\mathfrak{g}},\ldots,X_n^{\mathfrak{g}}$ of $\mathfrak{g}$, denote by $C^k_{ij}\in\genfd$ for $i,j,k\in\{ 1,\ldots,n\}$
the structure constants defined by
\begin{equation}\label{eq:Cijl}
  [ X_i^{\mathfrak{g}}, X_j^{\mathfrak{g}}] = \sum_{k = 1}^n C^k_{ij} X_k^{\mathfrak{g}}.
\end{equation}
Constants $C^k_{ij}$ are antisymmetric in lower indices and satisfy a quadratic relation reflecting the Jacobi identity in $\mathfrak{g}$. 
According to~\cite{ldWeyl}, there is a unique monomorphism of $\genfd$-algebras
$\iota: U({\mathfrak{g}})\to\hat{A}_{n,\genfd}$ extending the formulas 
\begin{equation}\label{eq:univfla}
X_i^{\mathfrak{g}} \mapsto \iota( X_i^{\mathfrak{g}}) = \sum_{l = 1}^n x_l \sum_{N = 0}^\infty \frac{(-1)^N}{N!} B_N (\bbC^N)^l_i,
\end{equation}
where $B_N$ is the $n$-th Bernoulli number and $\bbC$ is an $n\times n$ matrix
with values in $\genfd$, defined by
$$
\bbC^i_j = \sum_{k = 1}^n C^i_{jk} \partial^k.
$$
The monomorphism $\iota$ does not depend on the choice of the basis; over $\mathbb{R}$ and $\mathbb{C}$ the formula~(\ref{eq:univfla}) appeared to be known much before~(\cite{berezin1967,karasev}) and, suitably interpreted,
corresponds to the Gutt's star product~\cite{gutt}.
A simple differential geometric derivation of the formula~(\ref{eq:univfla})
over $\mathbb{R}$ is explained in detail in~\cite{halg}, Section~1.2. Similarly, Sections~7-9 of~\cite{ldWeyl} provide a geometrical derivation in formal geometry over any ring containing rationals. See also~\cite{Kathotia} for another point of view. Expression~(\ref{eq:univfla}) is related to the part of Campbell-Baker-Hausdorff series linear in the first argument~(\cite{ldWeyl}, Sections~7-9). Denote by
\begin{equation}\label{eq:symmap}
e^{\mathfrak{g}}:\genfd[x_1,\ldots,x_n]\to U({\mathfrak{g}}),\,\,\,\,\,\,\,\,\,
x_{\alpha_1}\cdots x_{\alpha_k}\mapsto
\frac{1}{k !}\sum_{\sigma \in \Sigma(k)}
X_{\alpha_{\sigma 1}}^{\mathfrak{g}}\cdots X_{\alpha_{\sigma k}}^{\mathfrak{g}}
\end{equation}
the standard symmetrization (or coexponential) map (of vector spaces), where the symmetric group on $k$ letters is denoted $\Sigma(k)$. Via monomorphism $\iota$, the expression on the right-hand side
of~(\ref{eq:symmap}) can be interpreted in $\hat{A}_{n,\genfd}$.
If we apply the resulting element of $\hat{A}_{n,\genfd}$ on $1$ using (the formal completion of) the Fock action, we recover back the left-hand side of~(\ref{eq:symmap}).
In other words~(\cite{ldWeyl,leib}),
\begin{equation}
((\iota\circ e^{\mathfrak{g}})(q))\triangleright 1  = q,\,\,\,\,\,\,q\in\genfd[x_1,\ldots,x_n],
\end{equation}
where $\triangleright$ denotes the Fock action by differential operators. 

In this paper, it is proven that already the tensorial form, 
\begin{equation}\label{eq:univfla2}
X_i\mapsto \tilde{X}_i := \sum_{l = 1}^n x_l\sum_{N = 0}^\infty A_N (\bbC^N)^l_i,
\end{equation}
of the universal formula~(\ref{eq:univfla}), with $A_1 = 1$, 
guarantees in characteristic $0$
that precisely the {\it symmetrically} ordered noncommutative expressions
$$\frac{1}{k!}\sum_{\sigma\in\Sigma(k)} X_{\alpha_{\sigma(1)}}\cdots X_{\alpha_{\sigma(k)}},$$
interpreted via the embedding~(\ref{eq:univfla2}), and after acting upon the vacuum, recover back the commutative product $x_{\alpha_1}\cdots x_{\alpha_k}$.
The coefficients $A_N$ in~(\ref{eq:univfla2}) may be arbitrary for $N>0$ and $A_1 = 1$, instead of the choice $A_N = \frac{(-1)^N}{N!} B_N$ for all $N$, and $X_i$ may be generators of an arbitrary finitely generated associative $\genfd$-algebra $U$, instead of the motivating choice $X_i = X_i^{\mathfrak{g}}\in U({\mathfrak{g}})$.

Even more generally, we may replace $A_N (\bbC^N)^l_i$ in~(\ref{eq:univfla2}) by any expression of the form $p^{N-1,l}_{ij}(\partial^1,\ldots,\partial^n)\partial^j$ provided that $p^{N-1,l}_{ij}=p^{N-1,l}_{ij}(\partial^1,\ldots,\partial^n)$ is a homogeneous polynomial of degree $(N-1)$ in $\partial^1,\ldots,\partial^n$, antisymmetric under interchange of $i$ and $j$. Note that the previous case involving $U(\gg)$ may be recovered by setting
$$p^{N-1,l}_{ij} = \frac{(-1)^N B_N}{N!}\sum_{s = 1}^n (\bbC^{N-1})_s^l C^s_{ij}.$$ 
We do not discuss when the correspondence~(\ref{eq:univfla2}) (or its generalization involving $p^{N-1,l}_{ij}$) extends to a homomorphism $U\to\hat{A}_{n,\genfd}$ of algebras (in physics literature also called a {\em realization} of $U$). If $U$ is tautologically defined as the subalgebra of $\hat{A}_{n,\genfd}$ generated by the expressions $\tilde{X}_i\in\hat{A}_{n,\genfd}$, we alert the reader that the corresponding PBW type theorem often fails and the dimension of the space of degree $k>1$ noncommutative polynomials in $\tilde{X}_i$ generically exceeds the dimension of the space of symmetric polynomials of degree $k$. 

In the rest of the article below,
$X_i$-s are defined as elements in $\hat{A}_{n,\genfd}$ from the start,
hence we proceed without a distinction between $X_i$ and $\tilde{X}_i$.

\section{Results}

\begin{theorem}\label{mainth}
  Assume $\genfd$ is a field of characteristic different from $2$. Let 
  \begin{equation}\label{eq:embgen}
X_i = x_i + \sum_{l = 1}^n x_l\sum_{N = 1}^\infty \sum_{j = 1}^n p^{N-1,l}_{ij}(\partial^1,\ldots,\partial^n)\partial^j,\,\,\,\,\,\,i = 1,\dots,n,
  \end{equation}
  be $n$ distinguished elements of $\hat{A}_{n,\genfd}$, where $p^{N-1,l}_{ij}(\partial^1,\ldots,\partial^n)$ are arbitrary homogeneous polynomials of degree $(N-1)$ in $\partial^1,\ldots,\partial^n$, antisymmetric in lower indices $i,j$. Let $\alpha:\{1,\ldots,k\}\to\{1,\ldots,n\}$ be any function.
  Then, in the index notation, $\alpha_i = \alpha(i)$, 
  \begin{equation}\label{eq:mainfla}
\sum_{\sigma\in\Sigma(k)} X_{\alpha_{\sigma(1)}}\cdots X_{\alpha_{\sigma(k)}} \triangleright 1
    = k!\,x_{\alpha_1}\cdots x_{\alpha_k}.
  \end{equation}
  \end{theorem}

\vskip .03in
{\it Proof.}
  We prove the theorem by induction on degree $k$.
  For $k = 1$ all terms with $N\geq 1$ vanish, because we apply
  at least one derivative to $1$.

For general $k$, we write the sum~(\ref{eq:mainfla})
over all permutations in $\Sigma(k)$ in a different way.
We use the fact that the set of permutations of $n$ elements $\Sigma(n)$
is in the bijection with the set of pairs $(i, \rho)$ where
$0 \leq i \leq k$ and $\rho\in \Sigma(k-1)$. This can be done in many ways,
but we use this concrete simple-minded bijection
$$
(i,\rho)\mapsto \sigma, \,\,\,\,\,\,\sigma(k) := \left\lbrace
\begin{array}{ll}
i, & k=1,\\
\rho(k-1), & k>1 \mbox{ and } \rho(k-1)< i,\\
\rho(k-1)+1, & k>1 \mbox{ and } \rho(k-1)\geq i.
\end{array}
\right.
$$
For example, $(3,(2,3,1,5,4))\mapsto (3,2,4,1,6,5)$.

Define a bijection
$\Theta_i : \{ 1,\ldots, k-1\}\to \{ 1,\ldots,i-1,i+1,\ldots, k\}$ by
$$
\Theta_i(j) := \left\lbrace
\begin{array}{ll}
j, &  j< i,\\
j+1, & j\geq i.
\end{array}
\right.
$$
Clearly now $\sigma(j+1) = \Theta_i(\rho(j))$ for $1\leq j<k$.

We may thus renumber the sum
$$\sum_{\sigma \in \Sigma(k)}X_{\alpha_{\sigma(2)}}\cdots X_{\alpha_{\sigma(k)}}$$
as the double sum
$$\sum_{i = 1}^k X_{\alpha(i)} \cdot
\sum_{\rho\in\Sigma(k-1)} X_{(\alpha\circ \Theta_i)(\rho(1))}
\cdots X_{(\alpha\circ \Theta_i)(\rho(k-1))}$$
By the assumption of induction,
$$
\sum_{\rho\in\Sigma(k-1)} X_{(\alpha\circ \Theta_i)(\rho(1))}
\cdots X_{(\alpha\circ \Theta_i)(\rho (k-1))}\triangleright 1
= (k-1)!\,x_{(\alpha\circ \Theta_i)(1)}
\cdots x_{(\alpha\circ \Theta_i)(k-1)}
$$
The function $\Theta_i$ takes all values between $1$ and $k$ except $i$
exactly once.

Therefore, the left-hand side of~(\ref{eq:mainfla}) may be rewritten as
\begin{equation}\label{eq:firststep}
(k-1)!\sum_{i = 1}^k
X_{\alpha(i)} \triangleright
(x_{\alpha(1)}\cdots x_{\alpha(i-1)} x_{\alpha(i+1)}\cdots x_{\alpha(k)}).
\end{equation}
Substituting the expression~(\ref{eq:embgen})
for $X_{\alpha(i)}$ in~(\ref{eq:firststep})
we immeditaly observe two summands.
Let $\delta$ be the Kronecker symbol. Then the first summand is 
$$
(k-1)! \sum_{i=1}^k  \sum_{r=1}^n  x_r \delta^r_{\alpha(i)}
\cdot
(x_{\alpha(1)}\cdots x_{\alpha(i-1)} x_{\alpha(i+1)}\cdots x_{\alpha(k)})
= k!\, x_{\alpha(1)}\cdots x_{\alpha(k)},
$$
yielding the desired right-hand side for~(\ref{eq:mainfla}).
Hence for the step of induction on $k$ it is sufficient to
show that the remaining summand
$$
(k-1)! \sum_{i=1}^k \sum_{N = 1}^\infty \sum_{l=1}^n x_l
\sum_{s=1}^n p_{\alpha(i)s}^{N-1,l} \partial^{s}
(x_{\alpha(1)}\cdots x_{\alpha(i-1)} x_{\alpha(i+1)}\cdots x_{\alpha(k)})
$$
vanishes. This follows if for any $N>0$ the contribution
\begin{equation}\label{eq:essterm}
\sum_{i=1}^k  \sum_{s=1}^n p_{\alpha(i)s}^{N-1,l}
\partial^{s} 
(x_{\alpha(1)}\cdots x_{\alpha(i-1)} x_{\alpha(i+1)}\cdots x_{\alpha(k)}) = 0.
\end{equation}
Let $s\in\{1,\ldots,n\}$ and $M(s) = \{ j\in\{1,\ldots,i-1, i+1,\ldots,k\} | s = \alpha(j)\}$. By elementary application of partial derivatives,
\begin{equation}\label{eq:jneqi}
  \partial^{s}(x_{\alpha(1)}\cdots x_{\alpha(i-1)} x_{\alpha(i+1)}\cdots x_{\alpha(k)}) =
  \sum_{j\in M(s)} \,\,\,\prod_{r\in \{1,\ldots,k\} \backslash\{i,j\}} x_{\alpha(r)}.
\end{equation}
In particular, the contributions from $s\notin\{\alpha(1),\ldots,\alpha(i-1),\alpha(i+1),\ldots,\alpha(k)\}$, that is for $M(s)=\emptyset$,
vanish and $\partial^{s} 
(x_{\alpha(1)}\cdots x_{\alpha(i-1)} x_{\alpha(i+1)}\cdots x_{\alpha(k)}) = 0$.

Thus, for fixed $i$, the overall sum over all $s\in\{1,\ldots,n\}$ 
becomes a new sum over all $j\in\{1,\ldots,i-1,i+1,\ldots,k\}$
and each $j\neq i$ appears precisely once, namely for $s = \alpha(j)$. 
For fixed pair $(i,s)$, notice that the summands
do not depend on $j\in M(s)$, but we do not use this fact. By antisymmetry, 
$p^{N-1,l}_{\alpha(i)\alpha(i)}=0$ if $\operatorname{char}\genfd\neq 2$, hence
we are free to add any terms multiplied by $p^{N-1,l}_{\alpha(i)\alpha(i)}$.
For fixed $i$, we conclude
$$
\sum_{s=1}^n p_{\alpha(i)s}^{N-1,l} \partial^{s}
(x_{\alpha(1)}\cdots x_{\alpha(i-1)} x_{\alpha(i+1)}\cdots x_{\alpha(k)})
= \sum_{j=1}^k p_{\alpha(i)\alpha(j)}^{N-1,l}
\prod_{r\in \{1,\ldots,k\} \backslash\{i,j\}} x_{\alpha(r)}.
$$
Regarding that $\prod_{r\in \{1,\ldots,k\} \backslash\{i,j\}} x_{\alpha(r)}$ is a symmetric tensor in $i,j$, and $p_{\alpha(i)\alpha(j)}^{N-1,l}$ is antisymmetric under exchange of $i$ and $j$, their contraction must be zero,
$$
 \sum_{i=1}^k\sum_{j=1}^k p_{\alpha(i)\alpha(j)}^{N-1,l}
\prod_{r\in \{1,\ldots,k\} \backslash\{i,j\}} x_{\alpha(r)} = 0.
$$ 
Therefore,~(\ref{eq:essterm}) follows,
and consequently the step of induction on $k$.\,\,\,\,\,\,\,\,\,\,\,$\Box$
\vskip .02in

The reader may want to understand the reindexing and cancellation arguments following formula~(\ref{eq:jneqi}) on an example where $\alpha$ is not injective. Suppose $n=3$, $k=4$, and $\alpha$ sends $1,2,3,4$ to $1,3,3,2$ respectively. Then $\sum_{i=1}^4\sum_{s=1}^4 p_{\alpha(i)s}^{N-1,l}\partial^s\left(\prod_{r\neq i} x_{\alpha(r)}\right)$ has contributions as follows: for $i=1$ one obtains $\sum_s p_{1 s}^{N-1,l}\partial^s(x_3 x_3 x_2)= p_{1 2}^{N-1,l} x_3 x_3+2p_{13}^{N-1,l} x_3 x_2$, for $i = 2$ and $i = 3$ equal contributions $\sum_s p_{3 s}^{N-1,l}\partial^s(x_1 x_3 x_2) = p_{31}^{N-1,l}x_3 x_2 + p_{32}^{N-1,l}x_1 x_3 + p_{33}^{N-1,l} x_1 x_2$, and for $i = 4$ one obtains $\sum_s p_{2s}^{N-1,l}\partial^s(x_1 x_3 x_3) = p_{21}^{N-1,l} x_3 x_3 + 2 p_{23}^{N-1,l} x_1 x_3$. By the antisymmetry of $p^{N-1,l}$, the double sum is $0$.  

\begin{corollary}
  Under the assumptions of Theorem~\ref{mainth},
  there is a well-defined $\genfd$-linear map
$$
  \tilde{e}:\genfd[x_1,\ldots,x_n]\to \hat{A}_{n,\genfd}
$$
  extending the formulas
  \begin{equation}\label{eq:Xsym}
    \tilde{e}: x_{\alpha_1}\cdots x_{\alpha_k}
    \mapsto \sum_{\sigma\in\Sigma(k)} X_{\alpha_{\sigma(1)}}\cdots X_{\alpha_{\sigma(k)}},
  \end{equation}
  for all $k\geq 0$ and for all (nonstrictly) monotone $\alpha:\{1,\ldots,k\}\to\{1,\ldots,n\}$. 
 Map $\tilde{e}$ satisfies
  \begin{equation}
  \tilde{e}(P_k)\triangleright 1 = k! P_k
  \end{equation}
  for all (commutative) polynomials $P_k = P_k(x_{\alpha_1},\ldots,x_{\alpha_n})$ homogeneous of degree $k$. In particular, $\tilde{e}$
  is injective iff $\operatorname{char}\genfd = 0$. In that case, the elements $e(x_{\alpha_1}\cdots x_{\alpha_n})$ are linearly independent.
  If $\operatorname{char}\genfd = 0$,
  a modified map $e:\genfd[x_1,\ldots,x_n]\to \hat{A}_{n,\genfd}$
  with normalization on $k$-homogeneous elements given by  
  \begin{equation}\label{eq:nXsym}
    e: x_{\alpha_1}\cdots x_{\alpha_k}
    \mapsto \frac{1}{k!}
    \sum_{\sigma\in\Sigma(k)} X_{\alpha_{\sigma(1)}}\cdots X_{\alpha_{\sigma(k)}},
  \end{equation}
  is an injection. 
\end{corollary}

The map $\tilde{e}$ is well-defined because the right-hand side in~(\ref{eq:Xsym}) is symmetric in $\alpha_1,\ldots,\alpha_k$.
Formula~(\ref{eq:mainfla}) can be restated as $\tilde{e}(-)\triangleright 1 = k!\,\id$.
Note that the expressions~(\ref{eq:Xsym}) do not span an associative subalgebra, but only a subspace $e(\genfd[x_1,\ldots,x_n])$ of the subalgebra $\genfd\langle X_1,\ldots,X_n\rangle$ of $\hat{A}_{n,\genfd}$ generated by $X_1,\ldots,X_n$, in general. Denote by $\pi:\genfd\langle X_1,\ldots,X_n\rangle\to\genfd[x_1,\ldots,x_n]$ the vector space projection given by the Fock action on the vacuum vector $1\in\genfd[x_1,\ldots,x_n]$, that is $\pi(P) = P\triangleright 1$, $P\in\genfd\langle X_1,\ldots,X_n\rangle$. If $\operatorname{char}\genfd = 0$, the map $e$ can be viewed as a $\genfd$-linear section of the projection map $\pi$. In particular, $e$ is an isomorphism onto its own image and $\operatorname{Ker}\,\pi\oplus\operatorname{Im}\,e = \genfd\langle X_1,\ldots,X_n\rangle$. 

\vskip .1in {\bf Acknowledgements.} I have proved the result in 2006 at IRB; the final writeup has been finished and submitted in the first days of my Autumn 2019 stay at IH\'ES. In the final stage, I have been partly supported by the Croatian Science Foundation under the Project ``New Geometries for Gravity and Spacetime'' (IP-2018-01-7615). 

\bibliographystyle{amsplain}

\end{document}